\documentclass[10pt]{article}
\usepackage{amsmath, mathtools, amssymb, braket, amsthm, graphicx, titlefoot, color}
\usepackage{bm}
\usepackage{setspace}
\usepackage{caption}
\usepackage{wrapfig}
\usepackage{authblk}
\usepackage{float}
\usepackage[square,numbers]{natbib}
\usepackage{geometry}
\newcommand*{\affaddr}[1]{#1} 
\newcommand*{\affmark}[1][*]{\textsuperscript{#1}}

\title{Limit Cycle Analysis of 3-D Nonlinear systems}
\author{%
	Souma Mazumdar\affmark[1],
	Premashis Kumar\affmark[2]
	 and 
	 Gautam Gangopadhyay\affmark[3]  \\
	\affaddr{\affmark[1]Department of Theoretical Sciences}\\
	\affaddr{\affmark[2,3]Department of Chemical, Biological and Macro-molecular Sciences}\\
	\affaddr{S. N. Bose National Centre for Basic Sciences \\ Block - JD, Sector - III, Salt Lake City, Kolkata - 700 106} \\
	\thanks{\affmark[1]Email: \texttt{souma.mazumdar@bose.res.in}}
	\thanks{\affmark[2]Email: \texttt{pkmanager007@gmail.com}}
	\thanks{\affmark[3]Email: \texttt{gautam@bose.res.in}}
}
 
\date{}

\begin{document}
	\maketitle 
	
	\begin{abstract}
		Considering Limit Cycles as one of the limits of Lienard equation, an analyis analogous to centre manifold analysis has been done for a $3-D$ nonlinear system exhibiting Limit Cycle. A rigorous study on radius of the Limit Cycle orbit has been done by considering $\lambda-\omega$ equations for the particular system and subsequently converting the system equations from cartesian to polar form. It has been shown through an analysis analogous to Centre Manifold Analysis and reduction of the system dynamics on a lower dimensional space, the Limit cycle radius undergoes an increment change. One example is provided to support the theoretical predictions.
      % \noindent{\bf Word Count 5600}
       
	\end{abstract}

	\section{Introduction}
	Limit cycles\cite{ye1986theory,christopher2007limit,strogatz2018nonlinear} are interesting features of stability that often arise in Nonlinear Dynamical Systems\cite{wiggins2003introduction,strogatz2018nonlinear}. For $2-D$ nonlinear dynamical systems it is often observed in the phase portraits the systems move in closed orbits around a fixed point. Thus Limit Cycles give valueable information about the stability and the asymptotic behaviour of a system. Limit cycles can be of both types: Attracting and Repelling depending on the behaviour of the system for $t$ approaching positive or negative $\infty$ respectively. Numerous examples of nonlinear systems with Limit Cycle behaviour are abound in nature like the Nonlinear Van Der Pol Oscillators\cite{kanamaru2007van,tavazoei2009more,barbosa2007analysis}. So studying systems which have Limit Cycles are subjects of general interest in Nonlinear Dynamics. In this paper we have made an attempt to do so. \par
	 The purpose of the paper is a rigorous study on the radius of the Limit Cycles, the overall asymptotic behaviour of the system when with the standard two degrees of freedom we couple a third degree of freedom and following an analysis similar to Centre manifold Analysis\cite{wiggins2003introduction,liu2000application,namachchivaya2003centre}, a reduction of the dynamics on a reduced space is obtained. Then following standard approach of Limit cycle analysis we calculate the effect on radius on the Limit Cycle on reduced space due the additional coupled direction.\par
     With the beginning of discussions Limit Cycles are shown to be one of the two limits of Lienard equation\cite{villari1987qualitative,feng2002explicit,burton1965generalized}, the other limit being a Centre. When the middle term of the Lienard equation is equated to $0$ evaluated at the fixed point $(0,0)$ Limit Cycle comes out as one of the limits when the evaluated value is less than $0$ and the Centre comes out as another limit when the evaluated value is $0$. This fine difference between a Limit cycle and a Centre in terms of their mathematical definitions as per the value of the middle term at the fixed point $(0,0)$ is brought out nicely through our work. Adding a third degree of freedom to the two dimensional Lienard equation and doing an analysis similar to the centre manifold analysis\cite{wiggins2003introduction} we show mathematically the two criterion for obtaining a Limit cycle and Centre and thereby establish the true motivation behind our work.\par 
     The second motivation behind our work comes from analysing the effect on the Limit Cycle radius on addition of a coupled third direction. It is known that Limit Cycle radius undergoes change on coupling with an oscillator which results in change in the Limit Cycle area. Further the oscillations can be brought down to zero by changing the coupling strengths where the coupling parameters can be described as controls. We in our work show that by addition of a stable third direction the Limit Cycle radius undergoes an increment change where the eigenvalue of the stable direction acts as the control parameter. The equations established points towards the fact that by changing the eigenvalue of the stable direction we can change the Limit Cycle radius by a desired amount.\par
     Although the main focus and motivation of our work involve finding theoretical aspects of generating and controlling
     the limit cycle oscillation by carrying out a new analytical method, it is capable of satisfying different requirements
     of various well-known and promising technological applications due to its general applicability.\par 
     This study can be implemented in designing a controller for manipulating the amplitude of the limit cycle oscillation 
     in various mechanical systems. The requirement of desired periodic motions in the manufacturing processes can be met 
     at low energy dissipation by adding control to the marginally stable limit cycle of a conservative 
     system\cite{schiehlen2005control}. The nonlinear behaviors of uncontrolled aeroelastic systems are well established. 
     Our finding in this report can be crucial in developing a reliable control technique of limit cycle 
     oscillations\cite{Strganac} in such aeroelastic systems. For example, our methodology can find direct application 
     in a two-dimensional aeroelastic system with nonlinear stiffness and revealing the role of stiffness on limit cycle 
     oscillation suppression of aeroelastic system by estimating the amplitude of oscillation more accurately under different 
     conditions. The amplitude of limit cycle oscillations in nonlinear aeroelastic systems\cite{shukla} can be minimized 
     by exploiting this method. The amplitude suppression of limit cycle oscillation can be achieved by a suitable choice 
     of control gains\cite{Nayfeh} and a time delay state feedback\cite{MACCARI2003123}. Here we have demonstrated the 
     possibility of having such phenomena in generic nonlinear systems by adding a coupled direction. Effective controllers 
     of the amplitude of the limit cycle in a general feedback scheme\cite{Tangamp} can be realized through our method of 
     considering the third degree of freedom coupled with the other two. Limit-cycle control and fuel control are popular
     for their usage in the automotive applications\cite{RIBBENS20131} and keeping the limit cycle oscillations' upper and 
     lower amplitude at a recommended value is at the heart of these controllers. Therefore we believe that our study of 
     controlling limit cycle oscillation can help in improving the performance of these controllers.\par 
     Besides these applications, many important biological processes like glycolysis\cite{selkov}, 
     circadian rhythm\cite{Goldbeter} follow limit cycle oscillation. Therefore studying peculiarities in limit cycle 
     oscillation would help to get more insights into those biological oscillations and understand the effect of different
     external and internal parameters on those processes. The general nature of our methodology makes it worthwhile enough 
     to study chemical oscillations from a different perspective. 
     In fact, the amplitude of the limit cycle plays a guiding role in the concentration dynamics of the system and thus 
     dictating the dissipation and energetics of the system. The general coupling term described in the report can be regarded
     as different coupling and controls that commonly exist in a model chemical system. Thus many amplitude-mediated emerging 
     and counterintuitive phenomena of chemical oscillatory systems in the presence of different coupling can be understood 
     at a fundamental level by implementing our approach.
     
     The paper is organised as follows. In the next section we bring up the Lienard equation and write it in a two dimensional form. Then adding a third direction we do an analysis similar to centre manifold analysis and show the two criterion mathematically for obtaining the Limit Cycle and a Centre. Then in subsequent sections we take up our analysis for a general $3-D$ system and on reduction on a two dimensional space was found to possess a Limit Cycle. Then bringing in the $\lambda-\omega$ \cite{murray2007mathematical,murray2001mathematical}  form for the reduced $2-D$ system  and converting the equations from cartesian to polar form we perform a deep investigation on the peculiarities of the radius of the Limit Cycle. As Limit Cycle and Centre are two limits of the Lienard Equation our motivation prompts us to do an analysis similar to centre manifold analysis which we term as "Limit Cycle Manfold Analysis". It is shown than the radius of the Limit Cycle orbit changes with an increment if the third added coupled direction is considered to be a stable direction. One example is provided to support the Theoretical predictions. Finally we conclude by adding some remarks.
     \section{Lienard Equation and its Limits}
     The starting point of our discussion is the Lienard Equation. The Lienard Equation\cite{burton1965generalized,villari1987qualitative,graef1972generalized} is a general equation which covers all the second order equations which have nonlinear damping termns.\\
     The standard form of Lienard equation is given by the following.
     \begin{equation}
     \begin{split}
     \ddot{x}+F(x,\dot{x})\dot{x}+G(x)=0
     \end{split}
     \end{equation}
     There are some prescribed conditions for a Lienard System to exhibit Limit Cycle. These are,
     \begin{itemize}
     	\item $xG(x) > 0$ for $|x|>0$
     	\item $\int_{0}^{\infty}G(x)dx=\int_{0}^{-\infty}G(x)dx=\infty$
     	\item $F(0,0)<0$
     	\item $\exists$ $x_{0}>0$ s.t. $F(x,\dot{x}) \geq 0$, $|x|\geq x_{0}$
     	\item $\exists$, $M>0$, s.t. $F(x,\dot{x})\geq-M$, $|x|\leq x_{0}$
     	\item $\exists$ $x_{1}>x_{0}$ s.t. $\int_{x_{0}}^{x_{1}}F(x,\dot{x})dx \geq 10Mx_{0}$ where $\dot{x}>0$ is an arbitrary decreasing positive function of $x$.
     \end{itemize}
     Among these we take the third condition that is $F(0,0)<0$ as the defining criteria for a Limit Cycle. \\
     Most of the second order nonlinear systems have $G(x)=x$. So for simplicity we have the value $x$ for $G(x)$. Writing as two first order systems Lienard Equation can be written as
     \begin{equation}
     \begin{split}
     &\dot{x}=y \\ &
     \dot{y}=-F(x,y)y-x
     \end{split}
     \end{equation}
     Now in the nonlinear term $F(x,y)$ we segregate the nonlinear part containing only $x,y$ and the constant part. We write $F(x,y) \rightarrow F(x,y)-k, (k>0)$. Then the Lienard Equation takes the form,
     \begin{equation}
     \begin{split}
     &\dot{x}=y \\ &
     \dot{y}=-[F(x,y)-k]y-x
     \end{split}
     \end{equation}
     On fixed point$(0,0)$ consideration $F(0,0)=0$. Futher if $k=0$ we get a Centre and if $k > 0$ we get a Limit Cycle.
     \subsection{Limit Cycle Manifold}
     Now we add a third direction coupled with the other two. The third direction is considered to be a stable direction. Then the system equations take the form,
     \begin{equation}
     \begin{split}
    & \dot{x}=y \\ &
    \dot{y}=ky-x-F(x,y)y \\ &
    \dot{z}=-\lambda z+f_{3}(x,y,z)
     \end{split}
     \end{equation}
     where $-\lambda(\lambda > 0)$ is the eigenvalue of the stable direction and $f_{3}(x,y,z)$ is the nonlinear term for the third direction. Now we add nonlinear couplings for the first and second direction and write the equations in matrix form as follows.
     \begin{equation}
     \begin{split}
     \begin{bmatrix}
     \dot{x} \\ \dot{y} \\ \dot{z}
     \end{bmatrix}=\begin{bmatrix}
     0 & 1 & 0 \\ -1 & k & 0 \\ 0 & 0 & -\lambda
     \end{bmatrix}\begin{bmatrix}
     x \\ y \\ z
     \end{bmatrix} + \begin{bmatrix}
     f_{1}(x,y,z) \\ f_{2}(x,y,z)-F(x,y) \\ f_{3}(x,y,z)
     \end{bmatrix}
     \end{split}
     \end{equation}
     where $f_{1}$ and $f_{2}$ are the nonlinear couplings for the first and second direction respectively. Note that the above system is in the block diagonal form for the linear part. Following an analysis similar to Centre Manifold analysis\cite{wiggins2003introduction} we consider $z$ as a polynomial of $x,y$. We write,
     \begin{equation}\label{Lim Man 1}
     \begin{split}
     z=h(x,y)=a_{0}x^2+a_{1}xy+a_{2}y^2+\mathcal{O}(3)
     \end{split}
     \end{equation}
     Differetiating with respect to $t$ we have the above equation as,
     \begin{equation}
     \begin{split}
     \dot{z}-\frac{\partial h}{\partial x}\dot{x}-\frac{\partial h}{\partial y}\dot{y}=0
     \end{split}
     \end{equation}
     Substituting the values of $\dot{x},\dot{y},\dot{z}$ from the matrix equations we have,
     \begin{equation}
     \begin{split}
     -\lambda z+f_{3}(x,y,z)-(2a_{0}x+a_{1}y)(y+f_{1}(x,y,z))-(2a_{2}y+a_{1}x)(-x+ky-F(x,y)+f_{2}(x,y,z))=0
     \end{split}
     \end{equation}
     Now for fixed point$(0,0)$ consideration we put $F(x,y)=0$ as it contains terms only in $x,y$. Then the above equation takes the form,
     \begin{equation}\label{mat eq expand 1}
     \begin{split}
     -\lambda z+f_{3}(x,y,z)-(2a_{0}x+a_{1}y)(y+f_{1}(x,y,z))-(2a_{2}y+a_{1}x)(-x+ky+f_{2}(x,y,z))=0
     \end{split}
     \end{equation} 
     Now we consider polynomial forms for $f_{1},f_{2},f_{3}$. We consider,
     \begin{equation}
     \begin{split}
    & f_{1}(x,y,z)=c_{0}x^2+c_{1}y^2+c_{2}z^2+c_{3}xy+c_{4}yz+c_{5}xz  \\ &
    f_{2}(x,y,z)=d_{0}x^2+d_{1}y^2+d_{2}z^2+d_{3}xy+d_{4}yz+d_{5}zx \\ &
    f_{3}(x,y,z)=e_{0}x^2+e_{1}y^2+e_{2}z^2+e_{3}xy+e_{4}yz+e_{5}zx
     \end{split}
     \end{equation}
     Substituting the values of $f_{1},f_{2},f_{3}$ in equation (\ref{mat eq expand 1}) we have,
     \begin{equation}
     \begin{split}
    & -\lambda z+(e_{0}x^2+e_{1}y^2+e_{2}z^2+e_{3}xy+e_{4}yz+e_{5}zx) \\ &-(2a_{0}x+a_{1}y)(y+c_{0}x^2+c_{1}y^2+c_{2}z^2+c_{3}xy+c_{4}yz+c_{5}zx)\\ &-(2a_{2}y+a_{1}x)(-x+ky+d_{0}x^2+d_{1}y^2+d_{2}z^2+d_{3}xy+d_{4}yz+d_{5}zx)=0
     \end{split}
     \end{equation}
     Substituting the value of $z$ from (\ref{Lim Man 1}) we have,
     \begin{equation}
     \begin{split}
     &-\lambda(a_{0}x^2+a_{1}xy+a_{2}y^2)+(e_{0}x^2+e_{1}y^2+e_{2}(a_{0}x^2+a_{1}xy+a_{2}y^2)^2+e_{3}xy+e_{4}y(a_{0}x^2+a_{1}xy+a_{2}y^2) \\ &+e_{5}x(a_{0}x^2+a_{1}xy+a_{2}y^2))-(2a_{0}x+a_{1}y)(y+c_{0}x^2+c_{1}y^2+c_{2}(a_{0}x^2+a_{1}xy+a_{2}y^2)^2+c_{3}xy \\ & +c_{4}y(a_{0}x^2+a_{1}xy+a_{2}y^2)+c_{5}x(a_{0}x^2+a_{1}xy+a_{2}y^2))-(2a_{2}y+a_{1}x)(-x+ky+d_{0}x^2+d_{1}y^2 \\ & +d_{2}(a_{0}x^2+a_{1}xy+a_{2}y^2)^2+d_{3}xy+d_{4}y(a_{0}x^2+a_{1}xy+a_{2}y^2)+d_{5}x(a_{0}x^2+a_{1}xy+a_{2}y^2))=0
     \end{split}
     \end{equation}
     Equating the coefficients of $x^2,y^2,xy$ to $0$ we from the above equation we have,
     \begin{equation}
     \begin{split}
     &-\lambda a_{0}+e_{0}+a_{1}=0 \\ &
     -\lambda a_{2}+e_{1}-a_{1}-2a_{2}k=0 \\ &
     -\lambda a_{1}+e_{3}-2a_{0}+2a_{2}-a_{1}k=0
     \end{split}
     \end{equation}
     Solving for $a_{0},a_{1},a_{2}$ we have,
     \begin{equation}
     \begin{split}
    & a_{0}=\frac{\lambda e_{0}+e_{3}+\frac{2(e_{0}+e_{1})}{\lambda+2k}+e_{0}k}{\lambda^{2}+2+\frac{2\lambda}{\lambda+2k}+\lambda k} \\ &
    a_{1}=\frac{\lambda e_{3}+\frac{2\lambda e_{1}}{\lambda+2k}-2e_{0}}{\lambda^2+2+\frac{2\lambda}{\lambda+2k}+\lambda k} 
    	\\ &
    	a_{2}=\frac{\lambda^2 e_{1}+2(e_{0}+e_{1})+\lambda ke_{1}-\lambda e_{3}}{(\lambda+2k)(\lambda^2+2+\frac{2\lambda}{\lambda+2k}+\lambda k)}
     \end{split}
     \end{equation}
     Choosing $e_{0}=e_{1}=0,e_{3}=1$ we have,
     \begin{equation}
     \begin{split}
     & a_{0}=\frac{1}{\lambda^2+2+\frac{2\lambda}{\lambda+2k}+\lambda k} \\ &
     a_{1}=\frac{\lambda}{\lambda^2+2+\frac{2\lambda}{\lambda+2k}+\lambda k} \\ &
     a_{2}=\frac{-\lambda}{(\lambda+2k)(\lambda^2+2+\frac{2\lambda}{\lambda+2k}+\lambda k)}
     \end{split}
     \end{equation}
     From equation (\ref{Lim Man 1}) we have,
     \begin{equation}\label{Manifold Eq}
     \begin{split}
     z&=a_{0}x^2+a_{1}xy+a_{2}y^2 \\ &
     =\frac{x^2}{\lambda^2+2+\frac{2\lambda}{\lambda+2k}+\lambda k}+\frac{\lambda xy}{\lambda^2+2+\frac{2\lambda}{\lambda+2k}+\lambda k}-\frac{\lambda y^2}{(\lambda+2k)(\lambda^2+2+\frac{2\lambda}{\lambda+2k}+\lambda k)}
     \end{split}
     \end{equation} 
     Puting $k=0$ in the above equation we get the Equation of Centre Manifold. Therefore the Centre Manifold Equation is,
     \begin{equation}
     \begin{split}
     z=\frac{x^2}{\lambda^2+4}+\frac{\lambda xy}{\lambda^2+4}-\frac{\lambda y^2}{\lambda(\lambda^2+4)}
     \end{split}
     \end{equation}
     Following the same spirit as in the case of Centre Manifold if we put $k > 0$ in equation in (\ref{Manifold Eq}) we get the Equation which we term as "Limit Cycle Manifold Equation". This forms the basis of our motivation to carry out Limit Cycle Analysis along the line of Centre Manifold Analysis.
     \section{Limit Cycle Manifold Analysis}
     As the purpose of our paper is to do the Limit Cycle Analysis of Dissipative systems we consider general systems which exhibit Limit Cycles which takes into account Limit Cycle behaviour of the Dissipative systems which is our point of interest. General systems which exhibit Limit Cycle behaviour can be cast into the following form\cite{murray2001mathematical,murray2007mathematical}.
     \begin{equation}
     \begin{split}
     & \dot{x}=\lambda(r)x-w(r)y \\ &
     \dot{y}=w(r)x+\lambda(r)y
     \end{split}
     \end{equation}
     where $\lambda(r),w(r)$ are functions of $r$ satisfying the conditions $\lambda(r) > 0$ when $r < r_{0}$ and $\lambda(r) < 0$ when $r > r_{0}$, $r_{0}$ being the Limit Cycle radius and $w(r) > 0$. It is usually chosen as,
     \begin{equation}
     \begin{split}
     &\lambda(r)=\gamma-r^2 \\ &
     w(r)=1
     \end{split}
     \end{equation}
     Let us try to find out the Limit cycle radius when the system is moving on the two dimensional plane. We have,
     \begin{equation}
     \begin{split}
      & \dot{x}=(\gamma-r^2)x-y \\ &
     \dot{y}=x+(\gamma-r^2)y
     \end{split}
     \end{equation}
     Therefore,
     \begin{equation}
     \begin{split}
    & x\dot{x}=(\gamma-r^2)x^2-xy \\ &
    y \dot{y}=xy+(\gamma-r^2)y^2
     \end{split}
     \end{equation}
     Adding the above two equations,
     \begin{equation}
     \begin{split}
     r\dot{r}=(\gamma-r^2)r^2
     \implies \dot{r}=(\gamma-r^2)r
     \end{split}
     \end{equation}
     Now the Limit Cycle condition is given by $\dot{r}=0$ which gives $r=\sqrt{\gamma}$. This is the Limit Cycle radius which we denote by $R_{0}$. \\
     Now adding a third direction which we considered as a stable direction as shown in the preceding section for Limit Cycle Manifold we have the system equations as,
     \begin{equation}
     \begin{split}
     & \dot{x}=(\gamma-(x^2+y^2))x-y+f_{1}(x,y,z) \\ &
     \dot{y}=x+(\gamma-(x^2+y^2))y+f_{2}(x,y,z) \\ &
     \dot{z}=-\lambda z+f_{3}(x,y,z)
     \end{split}
     \end{equation}
     where $f_{1},f_{2},f_{3}$, are the nonlinear parts of the above equations, $-\lambda(\lambda > 0)$ is the eigenvalue of the stable direction and $r^2=x^2+y^2$. Writing in matrix form the above equations take the form,
     \begin{equation}
     \begin{split}
     \begin{bmatrix}
     \dot{x} \\ \dot{y} \\ \dot{z}
     \end{bmatrix}=\begin{bmatrix}
     \gamma & -1 & 0 \\
      1 & \gamma & 0 \\
      0 & 0 & -\lambda
     \end{bmatrix}
     \begin{bmatrix}
     x \\ y \\ z
     \end{bmatrix}
     + \begin{bmatrix}
     f_{1}(x,y,z)-x(x^2+y^2) \\
     f_{2}(x,y,z)-y(x^2+y^2) \\ 
     f_{3}(x,y,z)
     \end{bmatrix}
     \end{split}
     \end{equation}
     Again linear part of the above system is in block diagonal form. The Linear part of the Limit Cycle equations(first two equations) is decoupled from the stable direction(third equation). Considering a "Limit Cycle Manifold" we write $z$ as a polynomial function of $x,y$.
     \begin{equation}\label{Limit Cycle Manifold}
     z=h(x,y)=a_{0}x^2+a_{1}xy+a_{2}y^2+\mathcal{O}(3)
     \end{equation}
     Differentiating the above equation with respect to $t$ we have,
     \begin{equation}
     \dot{z}-\frac{\partial h}{\partial x}\dot{x}-\frac{\partial h}{\partial y}\dot{y}=0
     \end{equation}
     Substituting the values of $\dot{x},\dot{y},\dot{z}$ in the above equation we get,
     \begin{equation}
     \begin{split}
     &-\lambda z+f_{3}(x,y,z)-(2a_{0}x+a_{1}y)[(\gamma-(x^2+y^2))x-y+f_{1}(x,y,z)]\\ &-(a_{1}x+2a_{2}y)[x+(\gamma-(x^2+y^2))y+f_{2}(x,y,z)]=0
     \end{split}
     \end{equation}
     Substituting the value of $z$ from equation (\ref{Limit Cycle Manifold}) we have,
     \begin{equation}\label{Limit Cycle Manifold 1}
     \begin{split}
     &-\lambda(a_{0}x^2+a_{1}xy+a_{2}y^2)+f_{3}(x,y,z)-(2a_{0}x+a_{1}y)[(\gamma-(x^2+y^2))x-y+f_{1}(x.y,z)]\\ &-(a_{1}x+2a_{2}y)[x+(\gamma-(x^2+y^2))y+f_{2}(x,y,z)]=0
     \end{split}
     \end{equation}
     Now we assume polynomial forms for $f_{1},f_{2},f_{3}$. We write,
     \begin{equation}
     \begin{split}
     &f_{1}(x,y,z)=c_{0}x^2+c_{1}y^2+c_{2}z^2=c_{3}xy+c_{4}yz+c_{5}zx \\ &
     f_{2}(x,y,z)=d_{0}x^2+d_{1}y62+d_{2}z^2+d_{3}xy+d_{4}yz+d_{5}zx \\ &
     f_{3}(x,y,z)=e_{0}x^2+e_{1}y^2+e_{2}z^2+e_{3}xy+e_{4}yz+e_{5}zx
     \end{split}
     \end{equation}
     Substituting the values of $f_{1},f_{2},f_{3}$ from equation (\ref{Limit Cycle Manifold 1}) we have the above equation as,
     \begin{equation}
     \begin{split}
&     -\lambda(a_{0}x^2+a_{1}xy+a_{2}y^2)+e_{0}x^2+e_{1}y^2+e_{2}(a_{0}x^2+a_{1}xy+a_{2}y^2)^2+e_{3}xy+e_{4}y(a_{0}x^2+a_{1}xy+a_{2}y^2) \\ & +e_{5}x(z_{0}x^2+a_{1}xy+a_{2}y^2)-(2a_{0}x+a_{1}y)[(\gamma-(x^2+y^2))-y+c_{0}x^2+c_{1}y^2+c_{2}(a_{0}x^2+a_{1}xy+a_{2}y^2)^2 \\ &+c_{3}xy+c_{4}y(a_{0}x^2+a_{1}xy+a_{2}y^2)+c_{5}x(a_{0}x^2+a_{1}xy+a_{2}y^2)]-(a_{1}x+2a_{2}y)[x+(\gamma-(x^2+y^2))y+d_{0}x^2 \\ &+d_{1}y^2+d_{1}(a_{0}x^2+a_{1}xy+a_{1}y^2)+d_{3}xy+d_{4}y(a_{0}x^2+a_{1}xy+a_{2}y^2)+d_{5}x(a_{0}x^2+a_{1}xy+a_{2}y^2)]=0
     \end{split}
     \end{equation}
     Equating the coefficients of the $x^2,y^2,xy$ to $0$ from the above equation we have,
     \begin{equation}
     \begin{split}
     &-\lambda a_{0}+e_{0}-2a_{0}\gamma-a_{1}=0 \\ &-\lambda a_{2}+e_{1}+a_{1}-2a_{2}\gamma=0 \\ &
     -\lambda a_{1}+e_{3}+2a_{0}-2a_{1}\gamma-2a_{2}=0
     \end{split}
     \end{equation}
     Solving for $a_{0},a_{1},a_{2}$ we have,
     \begin{equation}
     \begin{split}
     & a_{0}=\frac{e_{0}(\lambda+2\gamma)+\frac{2(e_{0}+e_{1})}{\lambda+2\gamma}-e_{0}}{\lambda(\lambda+2\gamma)+2\gamma(\lambda+2\gamma)+4} \\ &
     a_{1}=e_{0}-\frac{e_{0}(\lambda+2\gamma)^2+2(e_{0}+e_{1})-e_{3}(\lambda+2\gamma)}{\lambda(\lambda+2\gamma)+2\gamma(\lambda+2\gamma)+4} \\ &
     a_{2}=\frac{e_{0}+e_{1}}{\lambda+2\gamma}-\frac{e_{0}(\lambda+2\gamma)^2+2(e_{0}+e_{1})-e_{3}(\lambda+2\gamma)}{(\lambda+2\gamma)[\lambda(\lambda+2\gamma)+2\gamma(\lambda+2\gamma)+4]}
     \end{split}
     \end{equation}
     Choosing $e_{0}=e_{1}=0,e_{3}=1$ we have,
     \begin{equation}
     \begin{split}
    & a_{0}=-\frac{1}{\lambda(\lambda+2\gamma)+2\gamma(\lambda+2\gamma)+4} \\ &
    a_{1}=\frac{\lambda+2\gamma}{\lambda(\lambda+2\gamma)+2\gamma(\lambda+2\gamma)+4} \\ &
    a_{2}=\frac{1}{\lambda(\lambda+2\gamma)+2\gamma(\lambda+2\gamma)+4}
     \end{split}
     \end{equation}
     Then the "Limit Cycle Manifold" Equation is given by,
     \begin{equation}\label{Limit Cycle Manifold Eq}
     \begin{split}
     z=-\frac{x^2}{\lambda(\lambda+2\gamma)+2\gamma(\lambda+2\gamma)+4}+\frac{(\lambda+2\gamma)xy}{\lambda(\lambda+2\gamma)+2\gamma(\lambda+2\gamma)+4}+\frac{y^2}{\lambda(\lambda+2\gamma)+2\gamma(\lambda+2\gamma)+4}
     \end{split}
     \end{equation}
     The reduced equations on the "Limit Cycle Manifold" is given by
     \begin{equation}
     \begin{split}
     & \dot{x}=(\gamma-(x^2+y^2))x-y+y[-\frac{x^2}{\lambda(\lambda+2\gamma)+2\gamma(\lambda+2\gamma)+4}+\frac{(\lambda+2\gamma)xy}{\lambda(\lambda+2\gamma)+2\gamma(\lambda+2\gamma)+4}\\ &+\frac{y^2}{\lambda(\lambda+2\gamma)+2\gamma(\lambda+2\gamma)+4}] \\ &
     \dot{y}=x+(\gamma-(x^2+y^2))y+x[-\frac{x^2}{\lambda(\lambda+2\gamma)+2\gamma(\lambda+2\gamma)+4}+\frac{(\lambda+2\gamma)xy}{\lambda(\lambda+2\gamma)+2\gamma(\lambda+2\gamma)+4}\\ &+\frac{y^2}{\lambda(\lambda+2\gamma)+2\gamma(\lambda+2\gamma)+4}]
     \end{split}
     \end{equation}
     on substituting the value $z$ from equation (\ref{Limit Cycle Manifold Eq}). \\
     Writing $x=R\cos\theta,y=R\sin\theta$ we have $\dot{x}=\dot{R}\cos\theta-R\sin\theta\dot{\theta},\dot{y}=\dot{R}\sin\theta+R\cos\theta \dot{\theta}$. Puting the values of $\dot{x}$ and $\dot{y}$ in the above two equations we have,
     \begin{equation}
     \begin{split}
     & \dot{R}\cos\theta-R\sin\theta\dot{\theta}=(\gamma-R^2)R\cos\theta-R\sin\theta \\ & +R\sin\theta\left[-\frac{R^2\cos^2\theta}{\lambda(\lambda+2\gamma+2\gamma(\lambda+2\gamma)+4)}+\frac{(\lambda+2\gamma)R^2\cos\theta\sin\theta}{\lambda(\lambda+2\gamma)+2\gamma(\lambda+2\gamma)+4}+\frac{R^2\sin^2\theta}{\lambda(\lambda+2\gamma)+2\gamma(\lambda+2\gamma)+4}\right] \\ &
     \dot{R}\sin\theta+R\cos\theta\dot{\theta}=R\cos\theta+(\gamma-R^2)R\sin\theta\\ &+R\cos\theta\left[-\frac{R^2\cos^2\theta}{\lambda(\lambda+2\gamma)+2\gamma(\lambda+2\gamma)+4}+\frac{(\lambda+2\gamma)R^2\sin\theta\cos\theta}{\lambda(\lambda+2\gamma)+2\gamma(\lambda+2\gamma)+4}+\frac{R^2\sin^2\theta}{\lambda(\lambda+2\gamma)+2\gamma(\lambda+2\gamma)+4}\right]
     \end{split}
     \end{equation}
    Solving for $\dot{R}$ and $\dot{\theta}$ from the above equations we have
    \begin{equation}
    \begin{split}
    \dot{R}&=(\gamma-R^2)R+R\sin2\theta\left[\frac{-R^2\cos2\theta+(\lambda+2\gamma)\frac{R^2\sin2\theta}{2}}{\lambda(\lambda+2\gamma)+2\gamma(\lambda+2\gamma)+4}\right] \\ &
    =(\gamma-R^2)R+R^3\left[\frac{\frac{-\sin4\theta}{2}-\frac{(\lambda+2\gamma)}{4}\cos4\theta+\frac{\lambda+2\gamma}{4}}{2[\lambda(\lambda+2\gamma)+2\gamma(\lambda+2\gamma)+4]}\right]
    \end{split}
    \end{equation}
    \begin{equation}\label{thetadot}
    \begin{split}
    \dot{\theta}& =1+R^2\cos2\theta\left[\frac{-\cos^2\theta+\sin\theta\cos\theta(\lambda+2\gamma)+\sin^2\theta}{\lambda(\lambda+2\gamma)+2\gamma(\lambda+2\gamma)+4}\right] \\ &
    =1+R^2\cos2\theta\left[\frac{-\cos2\theta+\frac{\sin2\theta}{2}(\lambda+2\gamma)}{\lambda(\lambda+2\gamma)+2\gamma(\lambda+2\gamma)+4}\right]
    \end{split}
    \end{equation}
    We plot $R-t$ as shown in the Figure 1.
    \begin{figure}[htb]
         \begin{center}
       		\includegraphics[width=15cm,height=12cm]{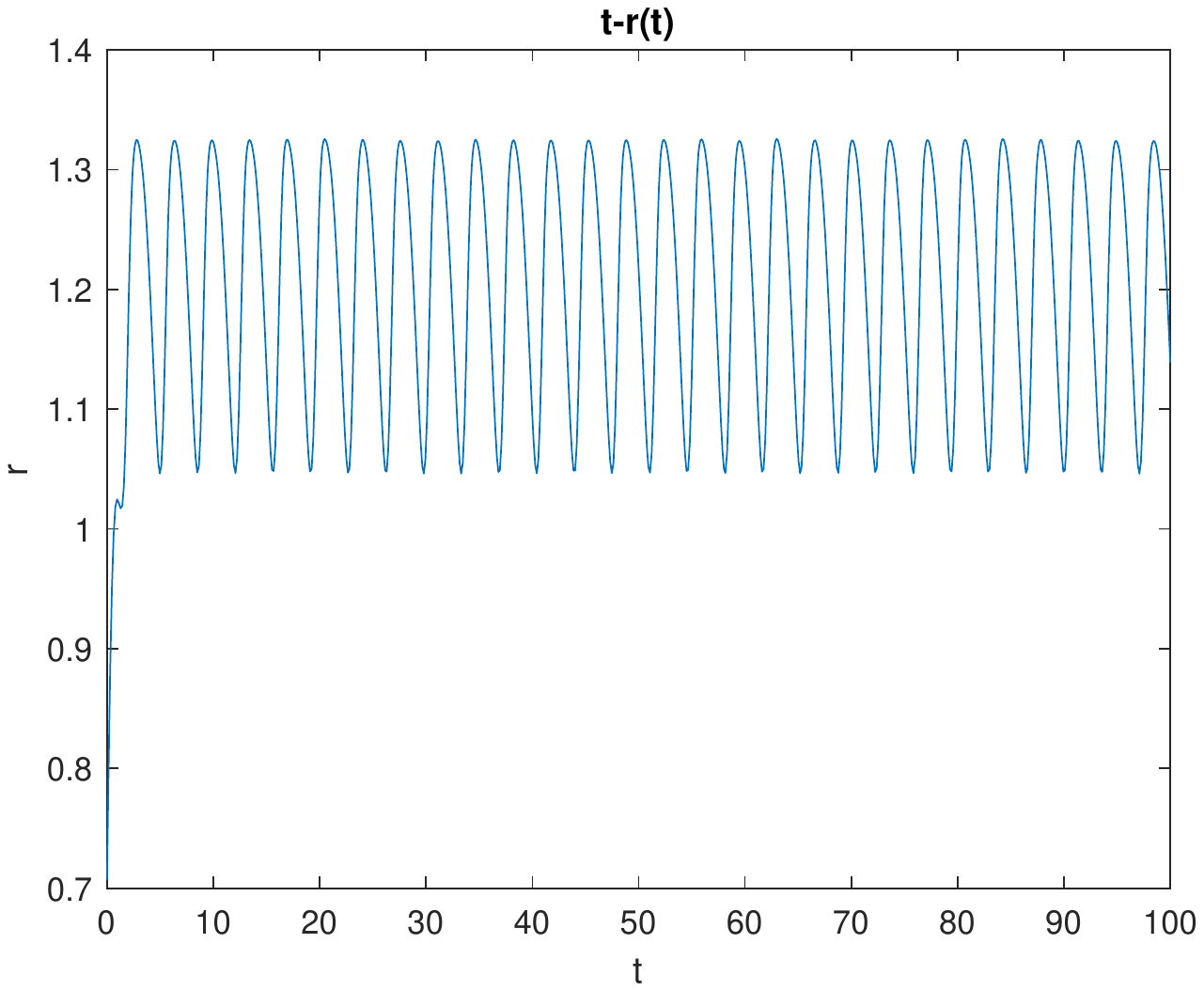}
      % 		\hspace{2mm}
       		\caption{\textbf{Limit Cycle Radius}}
    	\end{center}
       \end{figure}
      The $R-t$ plot shows an oscillatory motion about a mean radius which happens to be the new Limit cycle radius $R^{\prime}_{0}$. But the velocity is maximum on the mean radius$(R^{\prime}_{0})$. Thus we can write,
    \begin{equation}\label{pos VMAX}
    \begin{split}
    \dot{R}_{max}=\gamma R^{\prime}_{0}-R^{\prime3}_{0}+R^{\prime3}_{0}\left[\frac{-\frac{\sin4\theta}{2}-\frac{\lambda+2\gamma}{4}\cos4\theta+\frac{\lambda+2\gamma}{4}}{\lambda(\lambda+2\gamma)+2\gamma(\lambda+2\gamma)+4}\right]
    \end{split}
    \end{equation}
    Now we make an observation. We find that if we put $\theta+\frac{\pi}{2}$ in place of $\theta$, the motion is repeated. Therefore when $\theta$ equals $\theta+\frac{\pi}{4}$ the velocity should be maximum but in the opposite direction which is true for any oscillatory motion. So we can write,
    \begin{equation}\label{neg VMAX}
    \begin{split}
     -\dot{R}_{max}&=\gamma R^{\prime}_{0}-R^{\prime3}_{0}+R^{\prime3}_{0}\left[\frac{-\frac{\sin4(\theta+\frac{\pi}{4})}{2}-\frac{\lambda+2\gamma}{4}\cos4(\theta+\frac{\pi}{4})+\frac{\lambda+2\gamma}{4}}{\lambda(\lambda+2\gamma)+2\gamma(\lambda+2\gamma)+4}\right] \\ &
     =\gamma R^{\prime}_{0}-R^{\prime3}_{0}+R^{\prime3}_{0}\left[\frac{\frac{\sin4\theta}{2}+\frac{\lambda+2\gamma}{4}\cos4\theta+\frac{\lambda+2\gamma}{4}}{\lambda(\lambda+2\gamma)+2\gamma(\lambda+2\gamma)+4}\right]
    \end{split}
    \end{equation}
    Adding equations(\ref{pos VMAX}) and (\ref{neg VMAX}) we get,
    \begin{equation}
    \begin{split}
    \gamma R^{\prime}_{0}-R^{\prime3}_{0}+R^{\prime3}_{0}\frac{\frac{\lambda+2\gamma}{4}}{\lambda(\lambda+2\gamma)+2\gamma(\lambda+2\gamma)+4}=0
    \end{split}
    \end{equation}
    Therefore,
    \begin{equation}\label{Mean Increased Radius}
    \begin{split}
    R^{\prime}_{0}=\sqrt{\frac{\gamma}{1-\frac{\lambda+2\gamma}{4[\lambda(\lambda+2\gamma)+2\gamma(\lambda+2\gamma)+4]}}}
    \end{split}
    \end{equation}
    which gives the mean radius or the Limit Cycle radius on addition of a stable third direction. Without the third direction Limit Cycle radius is given by $R_{0}=\sqrt{\gamma}$. Thus we find that the Limit Cycle radius has undergone an increment change on addition of the stable direction. As a further check we see from equation (\ref{Limit Cycle Manifold Eq}) when the $z$ direction is absent that is $z=0$ which implies $\lambda(\lambda+2\gamma)+2\gamma(\lambda+2\gamma)+4=\infty$ we  have $R^{\prime}_{0}=R_{0}$. The increment of orbital radius on addition of a stable third direction as calculated above is a highly nontrivial result which our study reveals.
    \subsection{Asymptotic Amplitude and Angular Velocity} 
    \textbf{Asymptotic Amplitude:} \\ \\
    As the system settles in a limit cycle of mean radius $ R^{\prime}_{0}=\sqrt{\frac{\gamma}{1-\frac{\lambda+2\gamma}{4[\lambda(\lambda+2\gamma)+2\gamma(\lambda+2\gamma)+4]}}}$ the Asymptotic Amplitude is given by $R^{\prime}_{0}$. \\ \\
    \textbf{Asympotic Angular Velocity:} \\ \\
    The system moves in a mean circular path with phase increasing with time with a mean angular velocity which is the Asymptotic Angular Velocity. From equation (\ref{thetadot}) we have,
    \begin{equation*}
    \begin{split}
    \dot{\theta}
    =1+R^2\cos2\theta\left[\frac{-\cos2\theta+\frac{\sin2\theta}{2}(\lambda+2\gamma)}{\lambda(\lambda+2\gamma)+2\gamma(\lambda+2\gamma)+4}\right]
    \end{split}
    \end{equation*}
    Now the above equation is periodic with respect to $\theta$ with period $2\pi$. To get the mean angular velocity we have to take average on the L.H.S. and R.H.S. of the above equation\cite{sanders2007averaging}. Taking average on both sides of the above equation we have,
    \begin{equation}
    \begin{split}
    \dot{\bar{\theta}}&=\frac{1}{2\pi}\left[\int_{0}^{2\pi} d\theta-\frac{\bar{R}^2}{\lambda(\lambda+2\gamma)+2\gamma(\lambda+2\gamma)+4}\int_{0}^{2\pi}\cos^2 2\theta d\theta + \frac{\bar{R}^2}{2}\frac{\lambda+2\gamma}{\lambda(\lambda+2\gamma)+2\gamma(\lambda+2\gamma)+4}\int_{0}^{2\pi}\sin 2\theta \cos 2\theta d\theta\right] \\ &
    =\frac{1}{2\pi}\left[\int_{0}^{2\pi} d\theta- \frac{\bar{R}^2}{2[\lambda(\lambda+2\gamma)+2\gamma(\lambda+2\gamma)+4]}\int_{0}^{2\pi}(1+\cos 4\theta) d\theta+\frac{\bar{R}^2}{4}\frac{\lambda+2\gamma}{\lambda(\lambda+2\gamma)+2\gamma(\lambda+2\gamma)+4}\int_{0}^{2\pi}\sin 4\theta  d\theta \right] \\ &
    =\frac{1}{2\pi}\left[[\theta]_{0}^{2\pi}- \frac{\bar{R}^2}{2[\lambda(\lambda+2\gamma)+2\gamma(\lambda+2\gamma)+4]}\left[\theta+\frac{\sin 4\theta}{4}\right]_{0}^{2\pi}+\frac{\bar{R}^2}{4}\frac{\lambda+2\gamma}{\lambda(\lambda+2\gamma)+2\gamma(\lambda+2\gamma)+4}\left[-\frac{\cos 4\theta}{4}\right]_{0}^{2\pi} \right] \\ &
    =1-\frac{\bar{R}^2}{2[\lambda(\lambda+2\gamma)+2\gamma(\lambda+2\gamma)+4]}
    \end{split}
    \end{equation}
    Now $\bar{R}^2=R_{0}^{\prime 2}$. Puting the value of $R_{0}^{\prime 2}$ from equation (\ref{Mean Increased Radius}) in place of $\bar{R}^2$ we have,
    \begin{equation}
    \begin{split}
    &\dot{\bar{\theta}}=1-\frac{\frac{\gamma}{1-\frac{\lambda+2\gamma}{4[\lambda(\lambda+2\gamma)+2\gamma(\lambda+2\gamma)+4]}}}{2[\lambda(\lambda+2\gamma)+2\gamma(\lambda+2\gamma)+4]} \\ &
    =1-\frac{2\gamma}{4[\lambda(\lambda+2\gamma)+2\gamma(\lambda+2\gamma)+4]-(\lambda+2\gamma)}
    \end{split}
    \end{equation}
    \textbf{Time period for $1$ complete rotation around the Limit Cycle:} \\ \\
    Time period for $1$ complete rotation around the Limit Cycle is given by,
    \begin{equation}
    \begin{split}
    T&=\frac{2\pi}{\dot{\bar{\theta}}} \\ &
    \frac{2\pi}{1-\frac{2\gamma}{4[\lambda(\lambda+2\gamma)+2\gamma(\lambda+2\gamma)+4]-(\lambda+2\gamma)}}
    \end{split}
    \end{equation}
    \subsection{Total number of oscillations in one full cycle}
    We have shown velocity reverses to its maximum at phase separation of angle $\frac{\pi}{4}$. Therefore velocities are in the same phase when the phase separation is $\frac{\pi}{2}$. At this phase separation velocities are maximum in the same direction while crossing the mean radius. So we can say one oscillation is equivalent to phase difference $\frac{\pi}{2}$. Therefore the number of oscillations in $1$ complete cycle is given by,
    \begin{equation}
    \begin{split}
    \text{Number Of Oscillations in a full cycle}=\frac{2\pi}{\frac{\pi}{2}}=4
    \end{split}
    \end{equation}
    \subsection{Examples}
    We show two examples of Nonlinear System ocurring in nature which exhibit Limit Cycles. By the above analysis they will undergo a increase in Limit cycle radius on addition of a stable third direction.\\
    We cite the example of Van der Pol system which frequently occurs in nature.
    The $2-D$ model of this system is the following.\\ \\
    \textbf{Van der Pol:} \\
    \begin{equation}
    \begin{split}
   & \dot{x}=y  \\ &
    \dot{y}=-\mu(x^2-1)y-x
    \end{split}
    \end{equation}
    The O.D.E s in case of Van der Pol is given in polar form as,
   \begin{equation}
   \begin{split}
   &\dot{r}\cos\theta-r\sin\theta\dot{\theta}=r\sin\theta \\ &
   \dot{r}\sin\theta+r\cos\theta\dot{\theta}=-\mu(r^2\cos^2\theta-1)r\sin\theta-r\cos\theta
   \end{split}
   \end{equation}
   Solving for $\dot{r}$ we get,
   \begin{equation}\label{rdot}
   \begin{split}
   \dot{r}=-\mu r^3\cos^2\theta\sin^2\theta+\mu r\sin^2\theta
   \end{split}
   \end{equation}
   Now the above equation is periodic with respect to $\theta$ with period $2\pi$. To get the mean radius we average over the R.H.S and L.H.S of the above equation. Taking the average on both sides of equation (\ref{rdot}) we get,
   \begin{equation}
   \begin{split}
    \dot{\bar{r}}&=\frac{1}{2\pi}\left[-\mu \bar{r}^3\int_{0}^{2\pi}\cos^2\theta\sin^2\theta d\theta+\mu \bar{r}\int_{0}^{2\pi}\sin^2\theta d\theta\right] \\ &
   =\frac{1}{2\pi}\left[\frac{-\mu \bar{r}^3}{8}\int_{0}^{2\pi}(1-\cos 4\theta) d\theta+\frac{\mu \bar{r}}{2}\int_{0}^{2\pi}(1-\cos 2\theta) d\theta\right] \\ &
   =\frac{\mu \bar{r}}{2}\left[1-\frac{\bar{r}^2}{4}\right]
   \end{split}
   \end{equation}
   Therefore $\dot{\bar{r}}=0$ when $\bar{r}=2$. Therefore the mean radius of the Limit Cycle radius is given by $R_{0}=2$. Therefore $\gamma=4$. On addition of a stable third direction the Limit Cycle radius is given by
     \begin{equation}
    \begin{split}
    R^{\prime}_{0}=\sqrt{\frac{4}{1-\frac{\lambda+8}{4[\lambda(\lambda+8)+8(\lambda+8)+4]}}}
    \end{split}
    \end{equation}
    \textbf{Asymptotic Amplitude and Angular Velocity:} \\ \\
    The Asymptotic Amplitude is given by the above $R_{0}^{\prime}$. The Asymptotic Angular Velocity is given by
    \begin{equation}
    \begin{split}
    \dot{\bar{\theta}}=1-\frac{8}{4[\lambda(\lambda+8)+8(\lambda+8)+4]-(\lambda+8)}
    \end{split}
    \end{equation}
      We plot the $R-t$ plot for the Van der Pol systems as shown in Figure $2$. \\
    \begin{figure}[htb]
   	\centering
      	\caption{\textbf{Limit Cycle Radius: Van der Pol}}
   		\includegraphics[scale=1]{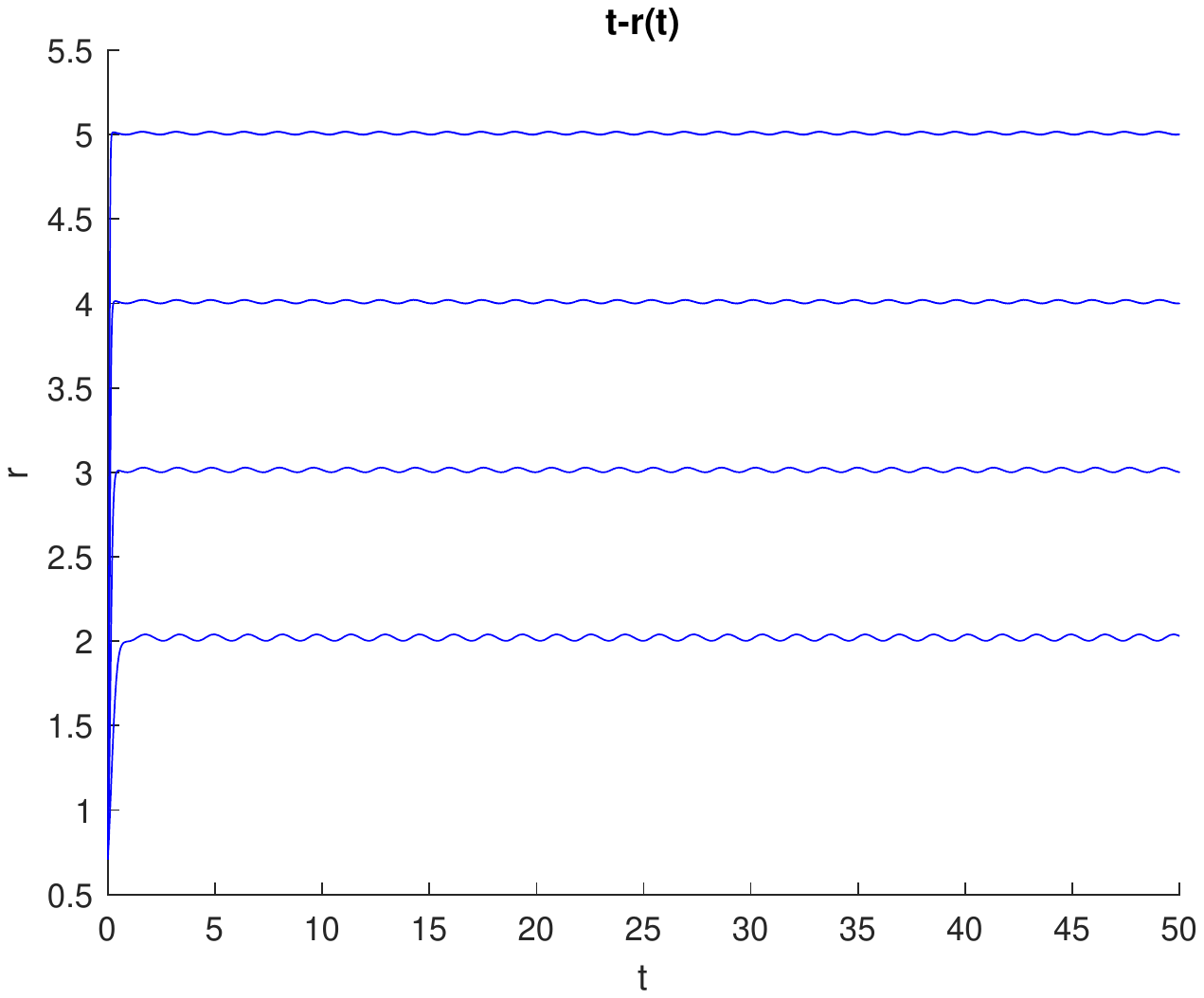}
     \end{figure}
    The figure shows the $R-t$ plots for Van der Pol system for different values of the parameter $\gamma$. The parameter values are chosen as $\gamma=4,9,16,25$. As $R_{0}=\sqrt{\gamma}$ the mean radius$(R_{0})=2,3,4,5$. On addition of a stable third direction, the mean radius$(R^{\prime}_{0})$ comes nearly equal to $R_{0}$ with difference in the second place of decimals. So the oscillations as shown in the plots is close to $R_{0}$. \\
    To see the dependence of modified radius with respect to parameter $\lambda$ we plot the modified mean radius $R^{\prime}_{0}$ with respect to $\lambda$ as shown in figure $3$. Here the value of the parameter $\gamma$ is kept fixed at $\gamma=4$.
    \begin{figure}[htb]
    	\centering
    	\caption{\textbf{Changing Limit Cycle radius with respect to parameters: Van der Pol}}
    	\includegraphics[scale=.75]{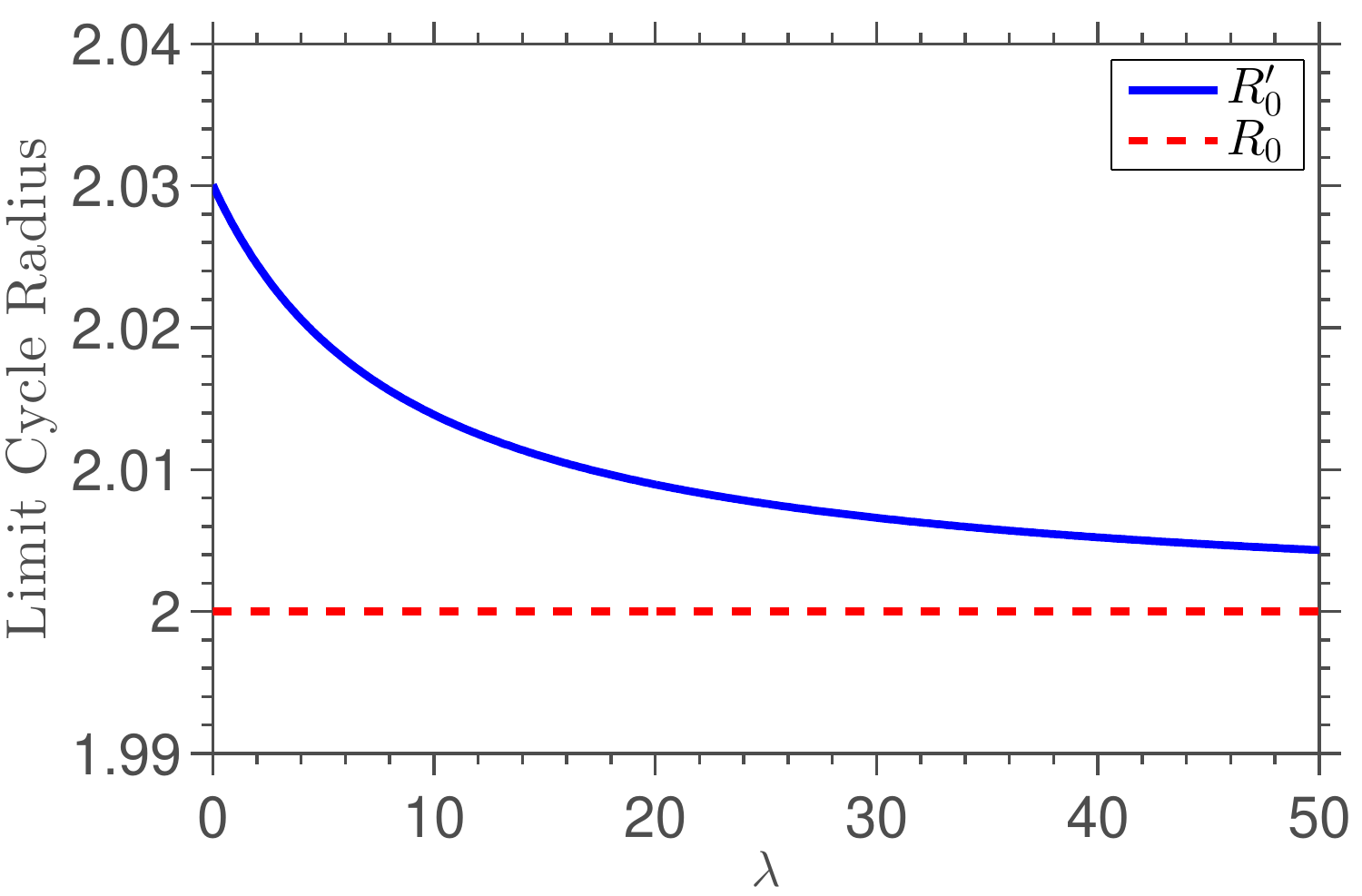}
    \end{figure}
  \section{Conclusion}
   We try to conclude the article with a short conclusion adding few points and remarks on our work. The first thing we like to point out is we have identified the Limit Cycle as one of the limits of Lienard Equation, the other limit being a Centre. This initiated the motivation behind carrying out an analysis following the centre manifold analysis which we have appropriately termed as "Limit Cycle Manifold" analysis. We have added a third direction  with our original two degrees of freedom and following an analysis analogous to centre manifold analysis reduced the system on a lower dimensional space. Through the initial discussion we successfully established the motivation behind our work. As our second motivation was to study the pecularities of Limit Cycle radius on addition of a coupled direction, we brought in the $\lambda-\omega$ system which is an inherent form for writing the differential equations for $2-D$ nonlinear systems exhibiting Limit Cycle. We did a detailed analysis on radius of the orbit of a Limit Cycle following the standard limit cycle analysis by converting the equations from cartesian to polar form. As we in our equation have added a third degree of freedom coupled with the other two we examined rigorously its effect on the Limit Cycle radius. What we have found that on addition of a stable direction the Limit Cycle radius changes with an increment which in our view is a completely new and interesting result that we have obtained theoretically. To tally with our theoretical prediction we have included three graphs which agree with the prediction.\par  
   
  \bibliographystyle{apa}
\bibliography{library}
\end{document}